\documentclass[a4paper,12pt,reqno,oneside]{amsart} % amsart loads amsmath, amsthm and amsfonts

\usepackage{setspace} 
\usepackage{typearea} % European margins
\usepackage{amscd,amssymb,verbatim} % amssymb loads amsfonts
\usepackage%[pdftex]
{graphicx} \setkeys{Gin}{width=\linewidth}
\usepackage{enumitem}
\usepackage{bbold}

\usepackage{xr-hyper}
\usepackage{cite} 
\usepackage[colorlinks,final,bookmarksnumbered,bookmarks]{hyperref}
\usepackage{amsrefs}

\usepackage{ifthen}
\newcommand{\arxiv}[2][]{\ifthenelse{\equal{#1}{}}
{\href{http://arxiv.org/abs/#2}{\tt arXiv:#2}}
{\href{http://arxiv.org/abs/math/#2}{\tt arXiv:math.#1/#2}}}

\usepackage[T1,T2A]{fontenc}
\usepackage[utf8]{inputenc}

\makeatletter
\renewcommand\subsection{\@startsection
{subsection}{2}{0cm} % name, level, indent
{-\baselineskip}     % before-skip
{0.5\baselineskip}   % after-skip
{\sffamily}} % style
\makeatother

\theoremstyle{plain}
\newtheorem{theorem}{Theorem}[section]

\newtheorem{corollary}[theorem]{Corollary}
\newtheorem{proposition}[theorem]{Proposition}

\theoremstyle{definition}

\newtheoremstyle{remark}
{}{}{}{}{\itshape}{}{ }{\thmname{#1}\thmnumber{ \itshape #2.}}
\theoremstyle{remark}

\newtheoremstyle{concise}
{}{}{}{}{\bfseries}{}{ }{\thmnumber{#2.}\thmnote{ #3.}}
\theoremstyle{concise}

\DeclareFontEncoding{LS1}{}{}
\DeclareFontSubstitution{LS1}{stix}{m}{n}
\DeclareFontEncoding{LS2}{}{}
\DeclareFontSubstitution{LS2}{stix}{m}{n}

\DeclareMathVersion{language}
\DeclareMathVersion{metalanguage}

\def\formulas{\mathversion{language}}

\def\metameta{\mathversion{normal}}

\DeclareMathAlphabet{\mf}{OML}{zplm}{m}{it} % pazo's \mathit
\DeclareMathAlphabet{\ts}{OT1}{cmss}{m}{sl} % cm's \mathsfit (no lowercase greek?)

\DeclareMathAlphabet{\fm}{U}{eur}{m}{n} % euler's roman (no slanted version)
\DeclareMathAlphabet{\tr}{OT1}{cmss}{m}{n} % standard (cm's) \mathsf (no lowercase greek?)
\DeclareMathAlphabet{\mm}{OML}{cmm}{m}{it} % usual \mathnormal, makes uppercase Greek letters slanted

\DeclareMathAlphabet{\script}{LS1}{stixscr}{m}{n}
\DeclareMathAlphabet{\mathbbb}{U}{bbold}{m}{n}

\SetSymbolFont{letters}{language}{U}{eur}{m}{n}
\SetSymbolFont{letters}{metalanguage}{OML}{zplm}{m}{it}

\newcommand{\newsym}[5]{\fontfamily{#2}\fontencoding{#1}\fontseries{#3}\fontshape{#4}\selectfont\char#5}
\makeatletter
\newcommand{\newmathsymbol}[6]{#1{\@Pimathsymbol{#2}{#3}{#4}{#5}{#6}}}
\def\@Pimathsymbol#1#2#3#4#5{\mathchoice
  {\@Pim@thsymbol{#1}{#2}{#3}{#4}{#5}\tf@size}
  {\@Pim@thsymbol{#1}{#2}{#3}{#4}{#5}\tf@size}
  {\@Pim@thsymbol{#1}{#2}{#3}{#4}{#5}\sf@size}
  {\@Pim@thsymbol{#1}{#2}{#3}{#4}{#5}\ssf@size}}
\def\@Pim@thsymbol#1#2#3#4#5#6{\mbox{\fontsize{#6}{#6}\newsym{#1}{#2}{#3}{#4}{#5}}}
\makeatother

\def\lgroup{\newmathsymbol{\mathopen}{LS2}{stixex}{m}{n}{"DC}}
\def\rgroup{\newmathsymbol{\mathclose}{LS2}{stixex}{m}{n}{"DD}}
\newcommand{\mq}[1]{\lgroup #1\rgroup\,}

\def\mc#1{\mq{\oldstylenums{#1}}}

\def\prin{\boldsymbol\cdot\hskip1.5pt}

\def\cltop{\top}
\def\clbot{\bot}
\def\ab{{\mathchoice
  {\mbox{\larger[1]$\times$}}
  {\mbox{\larger[1]$\times$}}
  {\mbox{\larger[-2]$\times$}}
  {\mbox{\larger[-4]$\times$}}
}}
\def\triv{\checkmark}

\def\To{\ \longrightarrow\ }
\def\tofrom{\leftrightarrow}
\def\Tofrom{\ \longleftrightarrow\ }
\def\imp{\Rightarrow}

\def\turnstile{\vdash}
\def\Turnstile{\vDash}

\DeclareMathSymbol{\impord}{\mathord}{symbols}{41}
\DeclareMathSymbol{\mand}{\mathbin}{operators}{`\&}

\DeclareMathSymbol{\oc}{\mathord}{operators}{`!}
\DeclareMathSymbol{\wn}{\mathord}{operators}{`?}

\DeclareMathSymbol{\alpha}     {\mathalpha}{letters}{"0B}
\DeclareMathSymbol{\beta}      {\mathalpha}{letters}{"0C}
\DeclareMathSymbol{\gamma}     {\mathalpha}{letters}{"0D}
\DeclareMathSymbol{\delta}     {\mathalpha}{letters}{"0E}
\DeclareMathSymbol{\epsilon}   {\mathalpha}{letters}{"0F}
\DeclareMathSymbol{\zeta}      {\mathalpha}{letters}{"10}
\DeclareMathSymbol{\eta}       {\mathalpha}{letters}{"11}
\DeclareMathSymbol{\theta}     {\mathalpha}{letters}{"12}
\DeclareMathSymbol{\iota}      {\mathalpha}{letters}{"13}
\DeclareMathSymbol{\kappa}     {\mathalpha}{letters}{"14}
\DeclareMathSymbol{\lambda}    {\mathalpha}{letters}{"15}
\DeclareMathSymbol{\mu}        {\mathalpha}{letters}{"16}
\DeclareMathSymbol{\nu}        {\mathalpha}{letters}{"17}
\DeclareMathSymbol{\xi}        {\mathalpha}{letters}{"18}
\DeclareMathSymbol{\pi}        {\mathalpha}{letters}{"19}
\DeclareMathSymbol{\rho}       {\mathalpha}{letters}{"1A}
\DeclareMathSymbol{\sigma}     {\mathalpha}{letters}{"1B}
\DeclareMathSymbol{\tau}       {\mathalpha}{letters}{"1C}
\DeclareMathSymbol{\upsilon}   {\mathalpha}{letters}{"1D}
\DeclareMathSymbol{\phi}       {\mathalpha}{letters}{"1E}
\DeclareMathSymbol{\chi}       {\mathalpha}{letters}{"1F}
\DeclareMathSymbol{\psi}       {\mathalpha}{letters}{"20}
\DeclareMathSymbol{\omega}     {\mathalpha}{letters}{"21}
\DeclareMathSymbol{\varepsilon}{\mathalpha}{letters}{"22}
\DeclareMathSymbol{\vartheta}  {\mathalpha}{letters}{"23}
\DeclareMathSymbol{\varpi}     {\mathalpha}{letters}{"24}
\DeclareMathSymbol{\varrho}    {\mathalpha}{letters}{"25}
\DeclareMathSymbol{\varsigma}  {\mathalpha}{letters}{"26}
\DeclareMathSymbol{\varphi}    {\mathalpha}{letters}{"27}

\DeclareMathSymbol{\bbomega}{\mathalpha}{letters}{"7F}

\def\Ds{\mathcal{D}}

\def\F{\mathcal{F}}

\def\G{\mathcal{G}}

\def\L{\script{L}}

\def\0{\mathbbb{0}}
\def\1{\mathbbb{1}}

\def\x{\times}
\def\but{\setminus}

\def\phi{\varphi}
\def\emptyset{\varnothing}

\renewcommand{\:}{\colon}
\DeclareMathOperator{\Int}{Int}

\DeclareMathOperator{\Hom}{Hom}

\DeclareMathOperator\suchthat{:}

\begin{document}

\title{A joint logic of problems and propositions}
\author{Sergey A. Melikhov}
%\date{\today}
\address{Steklov Mathematical Institute of Russian Academy of Sciences,
ul.\ Gubkina 8, Moscow, 119991 Russia}
\email{melikhov@mi-ras.ru}

\begin{abstract}
In a 1985 commentary to his collected works, Kolmogorov informed the reader that his 1932 paper 
{\it On the interpretation of intuitionistic logic} ``was written in hope that with time, the logic of 
solution of problems [i.e., intuitionistic logic] will become a permanent part of a [standard] 
course of logic.
A unified logical apparatus was intended to be created, which would deal
with objects of two types --- propositions and problems.''
We construct such a formal system as well as its predicate version, QHC, 
which is a conservative extension of both the intuitionistic predicate calculus QH 
and the classical predicate calculus QC.
The axioms of QHC are obtained as a result of a simultaneous formalization of two well-known alternative explanations 
of intiuitionistic logic: 1) Kolmogorov's problem interpretation (with familiar refinements by Heyting 
and Kreisel) and 2) the proof interpretation by Orlov and Heyting, as clarified and extended by G\"odel.
\end{abstract}

\maketitle
%\tableofcontents

The present note is a short communication on the author's preprints \cite{M1}, \cite{M2}, aiming at 
a more readable exposition of their main ideas and results.

\section{What is a logic, formally?}

According to modern textbooks, a derivation system for a first-order logic is supposed to consist of axiom schemata and 
inference rules, where both the former and the latter use meta-variables (rather than variables of the language) and 
may involve verbal side conditions such as ``provided that $x$ is not free in $P$'' or ``provided that $t$ is free for $x$ in $P(x)$''.
In the early textbooks of Hilbert--Ackermann, Hilbert--Bernays and P. S. Novikov, a derivation system involved only 
individual axioms (i.e.\ formulas, not schemata); but its inference rules anyway used meta-variables and involved 
verbal side conditions.

We will use a third type of formalism, where the entire derivation system is an object of a formal language
(called the language of the meta-logic).
This formalism is a simplified version of Paulson's formal meta-logic \cite{Pau1}, used in the {\tt Isabelle} proof assistant.
In this approach, principles and rules involve no meta-variables and no verbal side conditions; thus they are
entirely formal syntactic objects, {\it and so it makes sense to speak of their semantics}.
Meta-variables and side conditions do not entirely disappear: they are present in the meta-rules of inference, 
which are the same for all first-order logics.

Let us briefly describe our version of the formalism (cf.\ \cite{M1}, \cite{M0}).
Assume that we are given the language $\L$ of a first-order logic, including the notion of a {\it formula} of $\L$.
We need $\L$ to contain predicate variables (as in the textbooks by Church, Hilbert--Ackermann, Hilbert--Bernays 
and Novikov), and not just predicate constants as in some modern textbooks (which nevertheless tend to use 
propositional variables in their treatments of propositional logics).
The natural place for predicate constants is then not in $\L$, but in the language of a theory over the logic.

{\it Meta-formulas} are built inductively out of formulas
by using the {\it meta-conjunction} $\mand$, the {\it meta-implication} $\imp$ and the universal 
{\it meta-quantifiers} over arbitrary individual and predicate variables of $\L$.
Meta-quantification over all individual variables that occur freely in a given meta-formula $\F$ is abbreviated 
by $\mc{1}\F$; and over all predicate variables (of any arity) that occur freely in $\F$ by $\mc{2}\F$.

A {\it rule} $F_1,\dots, F_m/ G$, where $F_1,\dots,F_m$ and $G$ are formulas, is an abbreviation for 
the meta-formula $\mc{2}\big(\mc{1} F_1\,\mand\dots\mand\,\mc{1} F_m\imp \mc{1} G\big)$.
A {\it principle} $\prin G$, where $G$ is a formula, is an abbreviation for the meta-formula $\mc{2}\mc{1} G$.
Thus rules with no premisses can be identified with principles.

A {\it derivation system} $\Ds$ is a meta-conjunction of finitely many principles and rules,
which are called its {\it laws}%
\footnote{These play the role axiom schemata, but technically they are neither schemata nor individual axioms.}
and {\it inference rules}.
Another derivation system $\Ds'$ (based on the same language) is {\it equivalent} to $\Ds$ if
the meta-formula $(\Ds\imp\Ds')\mand(\Ds'\imp\Ds)$ is derivable using the meta-rules of inference.
(The latter are the introduction and elimination rules for $\mand$, $\imp$ and the meta-quantifiers, 
in the style of natural deduction, along with an $\alpha$-conversion rule.
See the details in \cite{M1} or \cite{M0} and examples of how the meta-rules work in \cite{M0}.)

A {\it logic} is an equivalence class of derivation systems.
If $L$ is the logic determined by $\Ds$, we denote by $\turnstile _L\F$ (or simply $\turnstile\F$
when $L$ is clear from context) the judgement that the meta-formula $\Ds\imp\F$ is derivable using 
the meta-rules of inference.
The judgment $\turnstile(\F_1\mand\dots\mand\F_m)\imp\G$ is also abbreviated by $\F_1,\dots,\F_m\turnstile\G$.
When $F_1,\dots,F_m$ and $G$ are formulas, ``$F_1,\dots,F_m\turnstile G$'' has the same meaning as in 
the textbooks by Church, Enderton, Kolmogorov--Dragalin and Troelstra--van Dalen.
On the other hand, ``$F_1,\dots,F_m\turnstile G$'' as defined in the textbooks by Schoenfield and Mendelson 
has the meaning of $\mc{1}F_1,\dots,\mc{1}F_m\turnstile \mc{1}G$.
The latter relation on tuples of formulas determines derivability of rules, and hence the logic.

Rules $R_1,\dots,R_n$ are said to {\it imply} a rule $Q$ if $R_1,\dots,R_n\turnstile Q$ holds.
\formulas
For example, in intuitionistic logic $\prin\fm{\neg\neg\alpha\to\alpha}$ implies $\prin\alpha\lor\neg\alpha$, 
using that $\turnstile\prin\neg\neg(\beta\lor\neg\beta)$.
However, $\neg\neg\alpha\to\alpha\not\turnstile\alpha\lor\neg\alpha$,
since the premiss becomes derivable upon substituting $\neg\beta$ for $\alpha$, but
the conclusion does not.
\metameta
(In the traditional formalism, it can be said that the schema $\alpha\lor\neg\alpha$ 
is not derivable from the schema $\neg\neg\alpha\to\alpha$; this indicates that the traditional
representation of principles by schemata is sometimes misleading.)

\section{What is intuitionistic logic, informally?}

When devising a new logic, one first has to engage in informal considerations, or else risk 
ending up with a ``wrong'' (i.e.\ unmotivated) logic.
Our logic QHC will emerge naturally from an analysis of informal semantics of intuitionistic logic.

\subsection{The problem interpretation} \label{problem interpretation}

The idea of Kolmogorov was that the objects of intuitionistic logic ``are in reality not theoretical 
propositions but rather problems''.
Here is a slightly clarified form of his explanation of intuitionistic logic \cite{Kol}, incorporating 
minor improvements due to Heyting (1934). 
See \cite{M0}*{\S\ref{int:BHK}, \S\ref{int:about-bhk}} for a more detailed discussion.

Let consider a class of specific problems (e.g.\ geometric construction problems) which may have parameters 
that run over a fixed domain $D$.
The meaning of problems that involve connectives or quantifiers is explained as follows:

\begin{itemize}
\item a solution of $\Gamma\land\Delta$ consists of a solution
of $\Gamma$ and a solution of $\Delta$;

\item a solution of $\Gamma\lor\Delta$ consists of an explicit choice
between $\Gamma$ and $\Delta$ along with a solution of the chosen
problem;

\item a solution of $\Gamma\to\Delta$ is a {\it reduction} of $\Delta$ to
$\Gamma$; that is, a general method of solving $\Delta$ on the basis
of any given solution of $\Gamma$;

\item the {\it absurdity} $\ab$ has no solutions, and
$\neg\Gamma$ is an abbreviation for $\Gamma\to\ab$;

\item a solution of $\exists x\, \Theta(x)$ is a solution of $\Theta(x_0)$
for some explicitly chosen $x_0\in D$;

\item a solution of $\forall x\, \Theta(x)$ is a general method
of solving $\Theta(x_0)$ for each $x_0\in D$.
\end{itemize}

A key element here is the notion of a {\it general method}, which Kolmogorov explains as follows.
If $\Gamma(\script X)$ is a problem depending on the parameter $\script X$ ``of any sort'', then 
``to present a general method of solving $\Gamma(\script X)$ for every particular value of $\script X$''
should be understood as ``to be able to solve $\Gamma(\script X_0)$ for every given specific value of 
$\script X_0$ of the variable $\script X$ by a finite sequence of steps, known in advance (i.e.\ before 
the choice of $\script X_0$)''.

As observed in the Troelstra--van Dalen textbook, the above six clauses alone do not suffice to discern
intuitionistic logic from classical.
Indeed, if $|\Gamma|$ denotes the set of solutions of the problem $\Gamma$, 
then these clauses guarantee that:
\begin{itemize}
\item $|\Gamma\land\Delta|$ is the product $|\Gamma|\x|\Delta|$;
\item $|\Gamma\lor\Delta|$ is the disjoint union $|\Gamma|\sqcup|\Delta|$;
\item there is a map $\script F\:|\Gamma\to\Delta|\to\Hom(|\Gamma|,|\Delta|)$ into the set of all maps;
\item $|\ab|=\emptyset$;
\item $|\exists x\,\Theta(x)|$ is the disjoint union $\bigsqcup_{d\in D} |\Theta(d)|$;
\item there is a map $\script G\:|\forall x\,\Theta(x)|\to\prod_{d\in D} |\Theta(d)|$ into the product.
\end{itemize}
If we force $\script F$ to be the identity map, we obtain classical logic.
Indeed, $|\Gamma\lor\neg\Gamma|=|\Gamma|\sqcup\Hom(|\Gamma|,\emptyset)$ is never empty; 
thus $\Gamma\lor\neg\Gamma$ has a solution for each problem $\Gamma$.

The missing ingredient here is the following interpretation of the second-order {\it meta-quantifier} --- which unfortunately 
is not discussed in any textbooks, but is almost explicit in Kolmogorov's paper \cite{Kol}.
If $\Phi=\Phi(\gamma_1,\dots,\gamma_n)$ is a closed formula containing no predicate variables
other than the ones listed, then the principle $\prin\Phi$ is interpreted by the problem
``Find a general method of solving the problem $\Phi(\Gamma_1,\dots,\Gamma_n)$ for arbitrary specific
problems $\Gamma_1,\dots,\Gamma_n$ of appropriate arities''.

For example, $\prin\fm{\gamma\lor\neg\gamma}$ is interpreted by the problem ``Find a general method that for an arbitrary specific 
problem $\Gamma$ either solves $\Gamma$ or derives a contradiction from the assumption that a solution of $\Gamma$ is given''.
We do not know such a general method, and {\it this is why intuitionistic logic differs from classical.}

Kolmogorov also offered an interpretation of rules, which after some clarification 
(in order to separate syntax from semantics) applies to arbitrary meta-formulas \cite{M0}.

\subsection{The proof interpretation}\label{proof interpretation}

A somewhat different interpretation of intuitionistic logic was proposed by Orlov (1928) 
and independently by Heyting (1930, 31).
For them the objects of intuitionistic logic are not problems, but the assertions of constructive 
{\it provability} of propositions.
See \cite{M1}*{\S\ref{g1:letters1}} for a more detailed discussion.

G\"odel was able to formalize this interpretation in his 1933 translation of propositional intuitionistic 
logic into S4 (a premature version of this translation appears already in the paper by Orlov).
Later, in his sketch of a 1938 lecture \cite{Goe1}, G\"odel attempted to axiomatize also 
the more precise notion of {\it proof} in the sense of Orlov and Heyting; he speaks of proofs 
``understood not in a particular system, but in the absolute sense (that is, one can make it evident)''.

In this sketch G\"odel considers a ternary relation ``$zBp,q$, that is, $z$ is a derivation of $q$ from $p$''.
But in fact he also uses a binary relation ``$aBq$'' which is presumably meant to abbreviate
$aB\cltop,q$.
G\"odel's axioms for $B$ are as follows (literally):
\smallskip

\begin{itemize}
\item ``$zB\phi(x,y)\To \phi(x,y)$'';
\item ``$uBv\To u'B(uBv)$'';
\item ``$zBp,q\mand uBq,r\To f(z,u)Bp,r$'';
\item ``if $q$ has been proved and $a$ is the proof, [then] $aBq$ is to be written down''.
\end{itemize}
\medskip

Instead of attempting to clarify the meaning of this in G\"odel's original terms,
let us consider a very similar but more clearly described logic.
Namely, let S4pr be the extension of classical predicate logic with the following 
additional elements of the language:

\begin{itemize}
\item an operator \,$\suchthat$\, associating to every formula $F$ and every term $t$
a formula $t\suchthat F$;
\item a unary function $'$ that associates to every term $t$ a term $t'$;
\item a binary function $[\cdot]$ that associates to every two terms $s,t$ a term $s[t]$;
\item an operator $*$ that associates to every formula $F$ a term $*_F$,
\end{itemize}

\noindent
and with the following additional laws and inference rules:%
\footnote{A closely related, but more complex logic was studied by Art\"emov \cite{Ar1}.}

\formulas
\begin{enumerate}
\item $\prin\tr t\suchthat p\To p$
\item $\prin\tr t\suchthat p\To\tr t'\suchthat (\tr t\suchthat p)$
\item $\prin\tr s\suchthat (p\to q)\To (\tr t\suchthat p\to\tr s[\tr t]\suchthat q)$
\smallskip
\item $\dfrac{p}{*_p\suchthat p}$
\end{enumerate}
\medskip

We will need

\begin{proposition} The following principles and rules are derivable in S4pr:
\begin{enumerate}
\item[(1$'$)] $\prin \neg (\tr t\suchthat\clbot)$
\item[(1$''$)] $\prin \exists \tr t\ \tr t\suchthat p\To p$
\smallskip
\item[(1$'''$)] $\dfrac{\tr t\suchthat p}p$
\smallskip
\item[(2$'$)] $\prin\tr t\suchthat p\To\tilde{\tr t}\suchthat (\exists\tr t\ \tr t\suchthat p)$
\end{enumerate}
\end{proposition}
\metameta

\subsection{The Proclus--Kreisel principle}

We will see below that from a certain formal viewpoint, the Kolmogorov--Heyting problem interpretation and 
the Orlov--Heyting proof interpretation have a lot in common, but do not reduce to each other.
In the literature, however, they are traditionally conflated (which is not entirely unreasonable as Heyting's 1931 
paper contains elements of both).
The result of this chronic conflation became known as the ``BHK interpretation'' of 
intuitionistic logic, initially called so after Brouwer, Heyting and Kreisel, and later after 
Brouwer, Heyting and Kolmogorov.

The principal contribution of Kreisel was, in terms of the problem interpretation, to observe that
for the six clauses to make sense, one must include the stipulation that {\it every solution of 
a problem $\Gamma$ must be supplied with a proof that it is indeed a solution of $\Gamma$}.
See \cite{M0}*{\S\ref{int:BHK}} for a more detailed discussion.
In fact, the same principle was emphasized already by the ancient Greeks, particularly Proclus, in 
the context of geometric construction problems.
(Martin-L\"of, on the contrary, emphasized that this principle fails in his intuitionistic type theory.)

\section{QHC}

In what follows QH and QC will stand for predicate intuitionistic and classical logics.
We are going to construct a new logic, called QHC, which fits Kolmogorov's description \cite{Kol3} 
that was quoted in the abstract.

\subsection{Language}

The language of QHC contains:
\begin{itemize}
\item logical constants:
\begin{itemize}
\item classical {\it truth} and {\it falsity}, denoted $\cltop,\clbot$, and
\item their intuitionistic counterparts {\it triviality} and {\it absurdity}, denoted $\triv,\ab$;
\end{itemize}
\item for each $n\ge 0$ countably many $n$-ary predicate variables of two sorts:
\formulas
\begin{itemize}
\item {\it problem variables}, denoted by fancy Greek letters ($\alpha,\beta,\dots$) and 
\item {\it proper predicate variables}, denoted by fancy Latin letters ($p,q,\dots$);
\end{itemize}
\item countably many individual variables, written with the sans-serif font: $\tr{x,y,\dots}$;
\item classical and intuitionistic connectives and quantifiers (see below);
\item two new operators $\wn$ and $\oc$.
\metameta
\end{itemize}

The usual math font will be used for meta-variables.

Each formula of QHC will have one of two sorts: it must be either a {\it c-formula} or an {\it i-formula}.
An {\it atomic c-formula} is an expression of the form $p(x_1,\dots,x_n)$, where $p$ is an $n$-ary proper 
predicate variable and $x_1,\dots,x_n$ are individual variables.
An {\it atomic i-formula} is an expression of the form $\alpha(x_1,\dots,x_n)$, where $\alpha$ is 
an $n$-ary problem variable and $x_1,\dots,x_n$ are individual variables.
A {\it c-formula} is either $\cltop$ or $\clbot$ or an atomic c-formula or an expression of 
the form $\wn\Phi$, where $\Phi$ is an i-formula; or an expression built inductively out of those of 
the previous four types using the classical connectives $\land$, $\lor$, $\to$, $\tofrom$, $\neg$ and the classical 
quantifiers $\exists$, $\forall$.
An {\it i-formula} is either $\triv$ or $\ab$ or an atomic i-formula or an expression of 
the form $\oc F$, where $F$ is a c-formula; or an expression built inductively out of those of 
the previous four types using the intuitionistic connectives $\land$, $\lor$, $\to$, $\tofrom$, $\neg$ and 
the intuitionistic quantifiers $\exists$, $\forall$.
We do not distinguish graphically between classical connectives/quantifiers and their intuitionistic
counterparts, but it is always easy to tell which is which because the sort of every subformula that 
has one of the first four types (from the definitions of a c-formula and of an i-formula)
is clear from its appearance.

%A {\it primitive c-formula} is either $\cltop$ or $\clbot$ or an atomic c-formula or an expression of the form $\wn\Phi$, where 
%$\Phi$ is an i-formula.
%A {\it primitive i-formula} is either $\triv$ or $\ab$ or an atomic i-formula or an expression of  the form $\oc F$, where $F$ is a c-formula.
%Arbitrary c-formulas (resp.\ arbitrary i-formulas) are built inductively out of primitive c-formulas (resp.\ primitive i-formulas) 
%using the classical (resp.\ intuitionistic) connectives  $\land,\lor,\to,\tofrom,\neg$ and the classical (resp.\ intuitionistic) 
%quantifiers $\exists,\forall$.
%We do not distinguish graphically between classical connectives/quantifiers and their intuitionistic counterparts, but it is always 
%easy to tell which is which because the sort of every primitive subformula is clear from its appearance.

In what follows, i-formulas will be denoted by unambiguously Greek capital letters, c-formulas by 
unambiguously Latin capital letters, and formulas of QHC of unspecified sort by capital letters 
that according to \TeX\ are both Latin and Greek (e.g.\ $A$, $B$, $E$).

\subsection{Informal semantics}

Closed c-formulas are interpreted by propositions, and closed i-formulas by problems.
In more detail:
\begin{itemize}
\item $n$-ary proper predicate variables are interpreted by $n$-ary predicates;
\item $n$-ary problem variables are interpreted by problems with $n$ parameters;
\item classical connectives and quantifiers are interpreted as usual, via truth tables;
\item intuitionistic connectives and quantifiers are interpreted as in \S\ref{problem interpretation};
\item if $\Phi$ is an i-formula interpreted by a problem $\Gamma$ with $n$ parameters, then 
$\wn\Phi$ is interpreted by the $n$-ary predicate ``$\Gamma$ has a solution''.
\item if $F$ is a c-formula interpreted by an $n$-ary predicate $P$, then $\oc F$ is interpreted 
by the problem ``Prove $P$'' with $n$ parameters.
\end{itemize}

If $F=F(\gamma_1,\dots,\gamma_m,p_1,\dots,p_n)$ is a closed c-formula containing no predicate variables
other than the ones listed, then the principle $\prin F$ is interpreted by the proposition
``$F(\Gamma_1,\dots,\Gamma_m,P_1,\dots,P_n)$ holds for arbitrary problems
$\Gamma_1,\dots,\Gamma_m$ and predicates $P_1,\dots,P_n$ of appropriate arities''.
The judgment $\turnstile\prin F$ can be interpreted by the judgment that this proposition is true.

If $\Phi=\Phi(\gamma_1,\dots,\gamma_m,p_1,\dots,p_n)$ is a closed i-formula containing no predicate variables
other than the ones listed, then the principle $\prin\Phi$ is interpreted by the problem
``Find a general method of solving $\Phi(\Gamma_1,\dots,\Gamma_m,P_1,\dots,P_n)$ for arbitrary 
problems $\Gamma_1,\dots,\Gamma_m$ and propositions $P_1,\dots,P_n$ of appropriate arities''.
The judgment $\turnstile\prin\Phi$ can be interpreted by the judgment that this problem has a solution.

The interpretation of rules and other meta-formulas is more involved; see \cite{M2}.

\subsection{Derivation system}\label{deductive}

Some laws and inference rules of the QHC calculus are immediate:

\begin{itemize}
\item[(0a)] All laws and inference rules of QC (see \cite{M0}*{\S\ref{int:logics}}).
\item[(0b)] All laws and inference rules of QH (see \cite{M0}*{\S\ref{int:logics}}).
\end{itemize}

Here formulas of QC and QH are identified respectively with c-formulas and i-formulas containing 
no occurrences of $\wn$ and $\oc$.
Let us note that by using meta-rules of inference (namely, the introduction and elimination rules for 
second-order meta-quantifiers) we can apply the classical laws and inference rules to arbitrary c-formulas 
and the intuitionistic laws and inference rules to arbitrary i-formulas.

Let us discuss the remaining part of the derivation system.

The proposition ``$\Gamma$ has a solution'', denoted $\wn\Gamma$, can be rephrased in the notation of
\S\ref{problem interpretation} as ``$|\Gamma|\neq\emptyset$''.
Then it follows from the six clauses of the problem interpretation (see \S\ref{problem interpretation})
that the following propositions must be true for any problems $\Gamma$, $\Delta$ and any parametric 
problem $\Theta$:

\begin{itemize}
\item $\wn(\Gamma\land\Delta)\Tofrom\wn\Gamma\land\wn\Delta$
\item $\wn(\Gamma\lor\Delta)\Tofrom\wn\Gamma\lor\wn\Delta$
\item $\wn(\Gamma\to\Delta)\To(\wn\Gamma\to\wn\Delta)$
\item $\neg\wn\ab$
\item $\wn\exists x\,\Theta(x)\Tofrom\exists x\,\wn\Theta(x)$
\item $\wn\forall x\,\Theta(x)\To\forall x\,\wn\Theta(x)$
\end{itemize}

These propositions correspond to the following principles in the language of QHC:

\formulas
\medskip
\begin{enumerate}
\item[(1a)] $\prin\wn(\gamma\land\delta)\Tofrom \wn\gamma\land\wn\delta$
\item[(1b)] $\prin\wn(\gamma\lor\delta)\Tofrom \wn\gamma\lor\wn\delta$
\item[(1c)] $\prin\wn(\gamma\to \delta)\To (\wn\gamma\to\wn\delta)$
\item[(1d)] $\neg\wn\ab$
\item[(1e)] $\prin\wn\exists\tr x\,\theta(\tr x)\Tofrom \exists\tr x\,\wn\theta(\tr x)$
\item[(1f)] $\prin\wn\forall\tr x\,\theta(\tr x)\To\forall\tr x\,\wn\theta(\tr x)$
\end{enumerate}
\medskip
\metameta
This will be our first group of laws of QHC (some of them will turn out to be redundant).

Informally, (1d) is saying that $\ab$ is not just the hardest problem (as guaranteed by the
explosion principle, $\fm{\prin\ab\to\gamma}$), but a problem that has no solutions whatsoever.
This is just one example of how the problem interpretation is not entirely captured in QH, and 
is more fully captured in QHC.

A consequence of the Proclus--Kreisel principle is that a solution of a problem $\Gamma$ yields a proof of
the existence of a solution of $\Gamma$.
This weaker principle is expressible in the language of QHC, and we add it to the derivation system:

\begin{enumerate}
\smallskip
\item[(1g)] $\fm{\prin\gamma\to \oc\wn\gamma}$.
\end{enumerate}
\smallskip

The final part of the derivation system comes from the proof interpretation.
The principles and rules (1$'$), (1$''$),  (1$'''$), (2$'$), (3) and (4) of \S\ref{proof interpretation},
which are derivable in S4pr, have the following direct analogues in the language of QHC:

\formulas
\smallskip
\begin{enumerate}
\item[(2a)] $\prin\neg\oc\clbot$
\item[(2b)] $\prin\wn\oc p\to p$
\medskip
\item[(2c)] $\dfrac{\oc p}{p}$
\medskip
\item[(2d)] $\prin\oc p\to \oc\wn\oc p$
\item[(2e)] $\prin\oc(p\to q)\To(\oc p\to \oc q)$
\medskip
\item[(2f)] $\dfrac{p}{\oc p}$
\end{enumerate}

In fact, (2d) becomes redundant, as it follows from (1g).
The interplay between the first and second parts of the derivation system goes much further than this.
In particular, it can be shown that the laws (1a), (1b), (1d), (1e) and (1f) are also redundant.

\begin{proposition}
The derivation system (0a)--(2f) is equivalent to the following one:
\begin{enumerate}[start=0]
\item All laws and inference rules of QC and of QH
\smallskip
\item $\dfrac{p}{\oc p}$
\smallskip
\item $\dfrac{\alpha}{\wn\alpha}$
\medskip
\item $\prin\wn\oc p\to p$
\item $\prin\alpha\to\oc\wn\alpha$
\item $\prin\oc(p\to q)\To(\oc p\to \oc q)$
\item $\prin\wn(\alpha\to \beta)\To (\wn\alpha\to\wn\beta)$
\item $\neg\oc\clbot$
\end{enumerate}
\end{proposition}

\subsection{Basic properties}

\begin{proposition} The following principles are derivable in QHC:
\begin{enumerate}
\item $\oc\clbot\tofrom\ab$
\item $\wn\ab\tofrom\clbot$
\item $\prin\oc p\land\oc q\Tofrom \oc(p\land q)$
\item $\prin\oc p\lor \oc q\To \oc(p\lor q)$
\item $\prin\forall\tr x\,\oc p(\tr x)\Tofrom \oc\forall \tr x\,p(\tr x)$
\item $\prin\exists\tr x\,\oc p(\tr x)\To \oc\exists \tr x\,p(\tr x)$
\item $\prin\neg\alpha\tofrom\oc\neg\wn\alpha$
\item $\prin(\oc\wn\alpha\to\oc\wn\beta)\Tofrom\oc(\wn\alpha\to\wn\beta)$
\end{enumerate}
\end{proposition}
\metameta

\begin{proposition}\label{Galois} For an i-formula $\Phi$ and a c-formula $F$,
$\turnstile\wn\Phi\to F$ if and only if $\turnstile\Phi\to\oc F$.
\end{proposition}

Thus $\wn$ and $\oc$ descend to a Galois connection between the Lindenbaum poset of the equivalence classes of
i-formulas and that of the c-formulas.

\formulas
\begin{proposition} Writing $\Box:=\wn\oc$, the following are derivable in QHC:
\begin{enumerate}
\item $\prin\Box p\to p$
\item $\prin\Box p\to\Box\Box p$
\item $\prin\Box(p\to q)\To(\Box p\to\Box q)$
\smallskip
\item $\dfrac{p}{\Box p}$
\end{enumerate}
\end{proposition}

The extension of QC by the operator $\Box$ satisfying (1)--(4) is the modal logic QS4.

\begin{proposition} Writing $\nabla:=\oc\wn$, the following are derivable in QHC:
\begin{enumerate}
\item $\prin\alpha\to\nabla\alpha$
\item $\prin\nabla\nabla\alpha\to\nabla \alpha$
\item $\prin\nabla(\alpha\to \beta)\To(\nabla \alpha\to\nabla \beta)$
\item $\nabla\ab\to\ab$
\end{enumerate}
\end{proposition}

The extension of QH by the operator $\nabla$ satisfying (1)--(4) will be denoted QH4.
This intuitionistic modal logic was studied by Fairtlough--Walton \cite{FW} and Aczel \cite{Ac}
and later by Arte\"mov--Protopopesku \cite{AP}.
If we drop (4), the resulting logic was studied by Curry \cite{Cu1}*{p.\ 120} as early as 1950, by
Goldblatt (1979) and many others.

\begin{proposition} The following are derivable in QHC:
\begin{enumerate}
\item $\prin\nabla\alpha\to\neg\neg \alpha$
\item $\prin\neg\nabla\alpha\tofrom\neg \alpha$
\item $\prin\neg\alpha\tofrom\nabla\neg \alpha$
\item $\prin\Box(p\land q)\Tofrom\Box p\land\Box q$
\item $\prin\nabla(\alpha\land \beta)\Tofrom\nabla \alpha\land\nabla \beta$
\end{enumerate}
\end{proposition}
\metameta

\section{Syntactic interpretations}\label{syntactic interpretations}

\subsection{The $\Box$-interpretation}
The classical provability translation of QH in QS4 (see \cite{M0}) extends to the following syntactic 
{\it $\Box$-interpretation} of QHC in QS4, denoted by $A\mapsto A_\Box$:
\begin{itemize}
\item atomic i-formulas are prefixed by $\Box$;
\item intuitionistic connectives and quantifiers become classical (also, $\ab$, $\triv$ 
become $\clbot$, $\cltop$) and are prefixed by $\Box$ (in fact, only $\to$ and $\forall$ really need to be prefixed);
\item $\wn$ is erased,  $\oc$ is replaced by $\Box$.
\end{itemize}

\begin{theorem} If $A_1,\dots,A_n\turnstile _{QHC} A$, then $(A_1)_\Box,\dots,(A_n)_\Box\turnstile _{QS4} A_\Box$.
\end{theorem}

\begin{corollary} QHC is a conservative extension of QS4; and therefore also of QC.
\end{corollary}

It is not hard to reformulate the $\Box$-interpretation as an extension of G\"odel's original form (see \cite{M0})
of the provability translation:
\begin{itemize}
\item $\ab$, $\triv$ become $\clbot$, $\cltop$;
\item intuitionistic connectives and quantifiers become classical and are postfixed by $\Box$ 
(in fact, only $\lor$, $\exists$ and $\to$ really need to be postfixed);
\item $\oc$ is erased, $\wn$ is replaced by $\Box$.
\end{itemize}

\subsection{The $\neg\neg$-interpretation}
The classical $\neg\neg$-translation of QC in QH (see \cite{M0}) extends to the following syntactic {\it $\neg\neg$-interpretation} 
of QHC in QH, denoted by $A\mapsto A_{\neg\neg}$:
\begin{itemize}
\item atomic c-formulas are prefixed by $\neg\neg$;
\item classical connectives and quantifiers become intuitionistic (also $\clbot$, $\cltop$ become $\ab$, $\triv$)
and are prefixed by $\neg\neg$ (in fact, only $\lor$ and $\exists$ really need to be prefixed);
\item $\oc$ is erased, $\wn$ is replaced by $\neg\neg$.
\end{itemize}

\begin{theorem} If $A_1,\dots,A_n\turnstile _{QHC} A$, then 
$(A_1)_{\neg\neg},\dots,(A_n)_{\neg\neg}\turnstile _{QH} A_{\neg\neg}$.
\end{theorem}

\begin{corollary} QHC is a conservative extension of QH.
\end{corollary}

The author's question whether QHC is also a conservative extension of QH4 was recently answered affirmatively 
by A. Onoprienko \cite{On1}, \cite{On2}.

It is not hard to reformulate the $\neg\neg$-interpretation as an extension of Kuroda's version (see \cite{M0})
of the classical $\neg\neg$-translation:
\begin{itemize}
\item if the entire formula is a c-formula, then it is prefixed by a $\neg\neg$;
\item $\clbot$, $\cltop$ become $\ab$, $\triv$;
\item classical connectives and quantifiers become intuitionistic and are postfixed by $\neg\neg$ 
(in fact, only $\forall$ really needs to be postfixed);
\item $\wn$ is erased, $\oc$ is replaced by $\neg\neg$.
\end{itemize}

\subsection{The $\nabla$-interpretation}
The $\Box$-interpretation of QHC in QS4 composed with the embedding of QS4 in QHC, $\Box\mapsto\wn\oc$,
can be ``improved'' so as to preserve the sorts of formulas.
Namely, if $A$ is a formula of QHC, let $A_\nabla$ be the formula of QHC obtained from $A$ by prefixing 
all atomic i-formulas and all intuitionistic connectives and quantifiers by $\nabla:=\oc\wn$. 
(In fact, of the connectives and quantifiers only $\lor$ and $\exists$ really need to be prefixed.)

\begin{theorem} If $A_1,\dots,A_n\turnstile A$, then $(A_1)_\nabla,\dots,(A_n)_\nabla\turnstile A_\nabla$.
Moreover, if $A$ is a formula of QH, then $\turnstile A_\nabla$ implies $\turnstile A$.
\end{theorem}

In particular, this yields a non-trivial embedding of QH in QHC (actually, in QH4), whose image is
a ``proof-irrelevant copy'' of the intuitionistic logic.
On the other hand, by using models we disprove both $A\turnstile A_\nabla$ and 
$A_\nabla\turnstile A$ for formulas of QH.

\subsection{The $\Diamond$-interpretartion}
The $\neg\neg$-interpretation of QHC in QH can also be ``improved'' so as to preserve the sorts of formulas.
Namely, if $A$ is a formula of QHC, let $A_\Diamond$ be the formula of QHC obtained from $A$ by prefixing all classical connectives
and quantifiers, all atomic c-formulas, and all $\wn$'s by $\Box\Diamond$, where $\Box=\wn\oc$ and $\Diamond=\neg\Box\neg$.
(In fact, $\land$ does not really need to be prefixed, and it would suffice to prefix $\to$ and $\forall$ by just $\Box$.)

\begin{theorem} If $A_1,\dots,A_n\turnstile A$, then $(A_1)_\Diamond,\dots,(A_n)_\Diamond\turnstile A_\Diamond$.
Moreover, if $A$ is a formula of QC, then $\turnstile A_\Diamond$ implies $\turnstile A$.
\end{theorem}

The resulting non-trivial embedding of QC in QHC is Fitting's embedding of QC in QS4 \cite{Fi}.
On the other hand, $A\turnstile A_\Diamond$ is disproved by considering models.

\subsection{Refined $\neg\neg$-translation}
The following reformulation of the $\Diamond$-interpretation, when restricted to QC, can be regarded as 
a refinement of Kolmogorov's original $\neg\neg$-translation:

\begin{itemize}
\item prefix every classical $\to$ by a $\wn$ and postfix it by an $\oc$;
\item prefix every classical $\forall$ and $\land$ by a $\wn$, and postfix it by $\neg\oc\neg$;
\item prefix every classical $\lor$ and $\exists$ by $\neg\wn\neg$, and postfix it by an $\oc$;
\item replace every $\nabla$ by a $\neg\neg$;
\item replace by $\neg\oc\neg$ every $\oc$ that is followed by an atomic c-formula;
\item if the formula starts with $\wn$, replace that $\wn$ by $\neg\wn\neg$;
\item if the entire formula is an atomic c-formula, prefix it by a $\Diamond$.
\end{itemize}

\formulas
Indeed, the latter interpretation has the effect of expressing all classical connectives and quantifiers
in terms of the intuitionistic ones along with $\wn$, $\oc$ and the classical $\neg$. 
For instance, $\exists\tr x\,p(\tr x)$ is interpreted as $\neg\wn\neg\exists\tr x\,\oc p(\tr x)$, i.e.\
``it is impossible to derive a contradiction from a construction of an $\tr x$ along with a proof of $p(\tr x)$''.
From the viewpoint of mathematical practice this makes better sense than Kolmogorov's original $\neg\neg\exists\tr x\,p(\tr x)$.
\metameta

\subsection{The $\neg\neg$-translation cannot be improved}
The following diagram is easily seen to commute:
\[
\begin{CD}
QH@>\text{$\nabla$-translation}>>QH4\\
@V\fm{\prin\alpha\lor\neg\alpha} VV@VV\nabla\mapsto\neg\neg V\\
QC@>\text{$\neg\neg$-translation}>>QH.
\end{CD}
\]
We show that this diagram is a pushout.
In other words, the classical $\neg\neg$-translation of QC into QH cannot be improved in the following sense:
it does not factor through any logic obtained by adding a set of laws to QH4 that are collectively strictly
weaker than $\fm{\prin\nabla\alpha\tofrom\neg\neg\alpha}$.

\section{Principles} \label{principles section}

The laws and inference rules of QHC emerged from the earliest and best known informal interpretations of 
intuitionistic logic, and so are hardly controversial (as long as we accept the language of QHC).
But could it be that some other equally indisputable principle or rule, expressible in 
the language of QHC, is not derivable in QHC and so is ``missing'' from the derivation system?

\formulas
\subsection{Hilbert's No Ignorabimus Principle}
The principle of decidability for i-formulas, $\prin\wn(\gamma\lor\neg\gamma)$, will be referred to as 
``Hilbert's No Ignorabimus Principle'', or the {\it H-principle} for brevity.
The relevant writings by Hilbert, where he stressed the principle of ``decidability of every mathematical
problem'', are discussed in detail in \cite{M2}.

\begin{proposition} Each of the following principles is equivalent to the H-principle:
\begin{enumerate}
\item $\prin\wn(\neg\neg\alpha\to\alpha)$;
\item $\prin\wn\neg\alpha\tofrom\neg\wn\alpha$;
\item $\prin\nabla\alpha\tofrom\neg\neg\alpha$;
\item $\prin\wn(\nabla\alpha\lor\neg\nabla\alpha)$;
\item $\prin\neg\neg\nabla\alpha\to\nabla\alpha$;
\item $\prin\wn(\neg\neg\nabla\alpha\to\nabla\alpha)$;
\item $\prin\wn(\alpha\to\beta)\Tofrom(\wn\alpha\to\wn\beta)$;
\item $\prin(\Box p\to\Box q)\Tofrom\Box(\Box p\to\Box q)$;
\item $\prin\Box p\tofrom\Diamond\Box p$;
\item $\prin\Diamond p\tofrom\Box\Diamond p$.
\end{enumerate}
\end{proposition}

\begin{corollary} QHC extended by the H-principle is a conservative extension of QS5.
\end{corollary}

\subsection{Kolmogorov's Stability Principle}

The principle $\prin\neg\oc\neg p\to\oc p$ will be referred to as ``Kolmogorov's Stablity Principle'', 
or the {\it K-principle} for brevity.
\metameta
This is related to an opinion expressed by Kolmogorov in a letter to Heyting which, as argued in 
\cite{M1}*{\S6.2.2} (see also \cite{M2}*{\S3.1.2}), may be interpreted as the assertion that in 
a constructive framework, every sentence must be either a problem or a {\it stable proposition}, that is,
a proposition $P$ such that the problem $\neg\oc\neg P\to\oc P$ has a solution.
(As discussed in \cite{M1}*{\S4.1}, stable propositions are closely related to stable problems, 
i.e.\ problems $\Gamma$ such that the problem $\neg\neg\Gamma\to\Gamma$ has a solution.)
\formulas

\begin{proposition} Each of the following principles and rules is equivalent to the K-principle:
\begin{enumerate}
\item $\prin\wn(\neg\oc\neg p\to\oc p)$;
\item $\prin\wn(\oc p\lor\oc\neg p)$;
\item $\prin\Box p\tofrom p$;
\item $\prin\oc\neg p\tofrom\neg\oc p$;
\item $\prin\oc(p\to q)\Tofrom(\oc p\to\oc q)$;
\item $\prin\Box(p\to q)\Tofrom(\Box p\to\Box q)$;
\item $\prin\Box(p\to q)\Tofrom\Box(\Box p\to\Box q)$;
\item $\prin\Box(p\lor q)\Tofrom(\Box p\lor\Box q)$;
\item $\prin\exists\tr x\, p(\tr x)\Tofrom \neg\wn\neg\exists\tr x\, \oc p(\tr x)$;
\item $\prin\forall\tr x\, p(\tr x)\Tofrom \wn\forall\tr x\, \neg\oc\neg p(\tr x)$;
\item $\dfrac{\Diamond p}{p}$;
\item $\dfrac{\neg\oc p}{\neg p}$.
\end{enumerate}
\end{proposition}

Akin to (3) is the following ``Proof Constructivity Principle'' ({\it PC-principle}): $\prin\nabla\alpha\to\alpha$.

\begin{proposition} (a) Each of the following principles is equivalent to the K-principle:
\begin{enumerate}
\item $\prin(p\to q)\Tofrom \wn(\oc p\to\oc q)$;
\item $\prin p\lor q\Tofrom\wn(\oc p\lor \oc q)$;
\item $\prin\exists\tr x\, p(\tr x)\Tofrom \wn\exists\tr x\, \oc p(\tr x)$;
\item $\prin\forall\tr x\, p(\tr x)\Tofrom \wn\forall\tr x\, \oc p(\tr x)$.
\end{enumerate}

\noindent
(b) Each of the following principles is equivalent to the PC-principle:
\begin{enumerate}
\item $\prin(\alpha\to\beta)\Tofrom \oc(\wn\alpha\to\wn\beta)$;
\item $\prin\alpha\lor\beta\Tofrom\oc(\wn\alpha\lor \wn\beta)$;
\item $\prin\exists\tr x\,\alpha(\tr x)\Tofrom \oc\exists\tr x\, \wn\alpha(\tr x)$;
\item $\prin\forall\tr x\, \alpha(\tr x)\Tofrom \oc\forall\tr x\, \wn\alpha(\tr x)$.
\end{enumerate}
\end{proposition}

\begin{theorem} (a) QHC extended by any collection of consequences of the K-principle is a conservative extension of QH and of QC.

(b) QHC extended by any collection of consequences of the PC-principle is a conservative extension of QH and of QS4.
\end{theorem}

\subsection{Further principles}

Next we consider the following principles, each of them a consequence of either the K-principle or the PC-principle:
\begin{itemize}
\item ``Reducibility Principle'' ({\it R-principle}):
$\prin(\nabla\alpha\to\nabla\beta)\To \nabla(\alpha\to\beta)$;
\item ``Disambiguation Principle'' ({\it D-principle}):
$\prin\nabla(\alpha\lor\beta)\To (\nabla\alpha\lor\nabla\beta)$;
\item {\it $\forall_\nabla$-principle}: 
$\prin\nabla\forall\tr x\, \alpha(\tr x)\Tofrom \forall\tr x\, \nabla \alpha(\tr x)$;
\item {\it $\exists_\Box$-principle}:
$\prin\Box\exists\tr x\, p(\tr x)\Tofrom \exists\tr x\, \Box p(\tr x)$;
\item {\it $\forall_\Box$-principle}:
$\prin\forall\tr x\,\Box p(\tr x)\Tofrom \Box\forall\tr x\, p(\tr x)$;
\item {\it $\exists_\nabla$-principle}:
$\prin\exists\tr x\,\nabla\alpha(\tr x)\Tofrom\nabla\exists\tr x\,\alpha(\tr x)$.
\end{itemize}

Let us also consider the {\it ``PC-rule''}: $\dfrac{\nabla\alpha}\alpha$, which is easily seen to be equivalent 
to $\dfrac{\wn\alpha}{\alpha}$.

\begin{proposition} The PC-principle is equivalent to the meta-conjunction of the PC-rule and the R-principle.
\end{proposition}

8 of the 9 principles introduced so far admit a parallel description:

\begin{proposition} \

(a) \parbox{\linewidth}{\begin{enumerate}
\item[] $\prin\wn\oc(p\lor q)\Tofrom \wn(\oc p\lor\oc q)$ is equivalent to the K-principle;
\item[] $\prin\oc\wn\alpha\lor\oc\wn\beta\Tofrom\oc(\wn\alpha\lor\wn\beta)$ is equivalent to the D-principle.
\end{enumerate}}
\smallskip

(b) 
\parbox{\linewidth}{\begin{enumerate}
\item[] $\prin\wn\oc(p\to q)\Tofrom\wn(\oc p\to\oc q)$ is equivalent to the K-principle;
\item[] $\prin(\wn\oc p\to\wn\oc q)\Tofrom\wn(\oc p\to\oc q)$ is equivalent to the H-principle;
\item[] $\prin\oc\wn(\alpha\to\beta)\Tofrom\oc(\wn\alpha\to\wn\beta)$ is equivalent to the R-principle.
\end{enumerate}}
\smallskip

(c) \parbox{\linewidth}{\begin{enumerate}
\item[] $\prin\wn\oc\exists\tr x\, p(\tr x)\Tofrom \wn\exists\tr x\, \oc p(\tr x)$ is equivalent to the $\exists_\Box$-principle;
\item[] $\prin\exists\tr x\, \oc\wn\alpha(\tr x)\Tofrom \oc\exists\tr x\, \wn\alpha(\tr x)$ is equivalent to the $\exists_\nabla$-principle.
\end{enumerate}}
\smallskip

(d) \parbox{\linewidth}{\begin{enumerate}
\item[] $\prin\oc\wn\forall\tr x\,\alpha(\tr x)\Tofrom\oc\forall\tr x\,\wn\alpha(\tr x)$ is equivalent to the $\forall_\nabla$-principle;
\item[] $\prin\forall\tr x\,\wn\oc p(\tr x)\Tofrom\wn\forall\tr x\,\oc p(\tr x)$ is equivalent to the $\forall_\Box$-principle.
\end{enumerate}}
\end{proposition}

\begin{theorem} (a) The 9 previously introduced principles and the PC-rule (denoted PC$^*$) are not derivable in QHC.
Moreover, there are no implications between them, except for the self-implications and the following implications:
\smallskip

\begin{center}
\includegraphics[width=8cm]{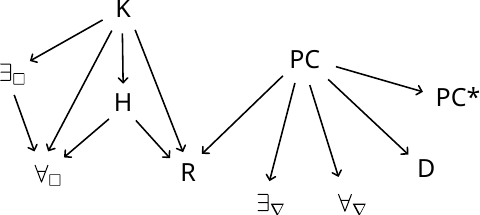}
\end{center}

(b) The following ``Exclusive Disambiguation Principle'' (ED-principle): 
\[\fm{\prin\neg(\alpha\land\beta)\To\big(\nabla(\alpha\lor\beta)\to\nabla\alpha\lor\nabla\beta\big)}\]
is not derivable in QHC, and moreover is not implied by any of the 9 principles, except the D-principle 
and the PC-principle.
\end{theorem}

Whether the PC-rule implies the ED-principle was posed as a problem by the author; it was recently answered 
in the negative by A. Onoprienko \cite{On1}.

\subsection{Intuitionistically unacceptable principles}

\formulas
\begin{proposition} Assume the H-principle. Then:

(a) the PC-principle is equivalent to the PC-rule and also to $\prin\alpha\lor\neg\alpha$;

(b) the D-principle is equivalent to Jankov's Principle, $\prin\neg\alpha\lor\neg\neg\alpha$;

(c) the $\forall_\nabla$-principle is equivalent to the Double Negation Shift Principle,
\[\prin\forall \tr x\,\neg\neg\alpha(\tr x)\To\neg\neg\forall \tr x\,\alpha(\tr x);\]

(d) the $\exists_\nabla$-principle is equivalent to the Strong Markov Principle,
\[\prin\neg\neg\exists \tr x\,\alpha(\tr x)\To\exists \tr x\,\neg\neg\alpha(\tr x);\]

(e) the ED-principle is equivalent to $\prin\nabla\alpha\lor\neg\nabla\alpha$.
\end{proposition}

\begin{proposition} Assume the PC-principle. Then:

(a) the $\forall_\Box$-Principle implies the Constant Domain Principle,
\[\prin\forall \tr x\,\big(\alpha\lor\beta(\tr x)\big)\To\alpha\lor\forall \tr x\,\beta(\tr x);\]

(b) the $\exists_\Box$-Principle implies the Principle of Independence of Premise,
\[\prin\big(\alpha\to\exists \tr x\,\beta(\tr x)\big)\To\exists \tr x\,\big(\alpha\to\beta(\tr x)\big).\]
\end{proposition}

\begin{proposition}\label{coupled principles}
(a) The principle $\prin\oc\exists\tr x\, p(\tr x)\Tofrom\exists\tr x\,\oc p(\tr x)$
\begin{itemize}
\item is equivalent to the meta-conjunction of the $\exists_\Box$- and $\exists_\nabla$-principles; and
\item implies the Generalized Markov Principle, 
$\prin\neg\forall \tr x\,\alpha(\tr x)\To\exists \tr x\,\neg \alpha(\tr x)$.
\end{itemize}

(b) The principle $\prin\wn\forall\tr x\,\alpha(\tr x)\Tofrom\forall\tr x\,\wn\alpha(\tr x)$
\begin{itemize}
\item is equivalent to the meta-conjunction of the $\forall_\Box$- and $\forall_\nabla$-principles; and
\item implies the Parametric Distributivity Principle,
\[\prin\neg\forall\tr y\,(\alpha(\tr y)\lor\forall\tr x\,\beta(\tr x,\tr y))\To
\neg\forall\tr y\,\forall\tr x\,(\alpha(\tr y)\lor\beta(\tr x,\tr y)).\]
\end{itemize}
\end{proposition}

\begin{proposition} Jankov's Principle, $\prin\neg\alpha\lor\neg\neg\alpha$, implies the ED-principle.
\end{proposition}
\metameta

It is well-known that Jankov's principle (also known as de Morgan's law or
``the weak law of excluded middle''), the Double Negation Shift Principle, the Strong Markov Principle,
the Generalized Markov Principle, the Constant Domain Principle and the Principle of Independence of Premise
are not derivable in QH (see \cite{M0}).
The Parametric Distributivity Principle is a variation of a well-known principle of Kleene; it is also not
derivable in QH \cite{M0}.

\section{Topological models}\label{models}

\subsection{Interior-based models} \label{interior-based}
By composing the $\Box$-interpretation of QHC with a topological model of QS4, we obtain 
an {\it interior-based model} of QHC (previously called an ``Euler--Tarski model'' 
in early versions of \cite{M1} and \cite{M2}).
Thus 
\begin{itemize}
\item we fix a topological space $X$ and a set $D$;
\item c-formulas with $n$ free variables are interpreted by $D^n$-indexed families of 
arbitrary subsets of $X$, as in the usual Leibniz--Euler--Venn models;
\item i-formulas with $n$ free variables are interpreted by $D^n$-indexed families of
open subsets of $X$, as in the usual Stone--Tang--Tarski models;
\item $\oc$ is interpreted by taking the interior, and $\wn$ by doing nothing.
\end{itemize}

\subsection{Regularization-based models} \label{regularization-based}
By composing the $\neg\neg$-interpretation of QHC with a Stone--Tang--Tarski model of QH, we get 
a {\it regularization-based model} of QHC (previously called a ``Tarski--Kolmogorov model'' 
in early versions of \cite{M1} and \cite{M2}).
Thus
\begin{itemize}
\item we fix a topological space $X$ and a set $D$;
\item i-formulas with $n$ free variables are interpreted by $D^n$-indexed families of
open subsets of $X$, as in the usual Stone--Tang--Tarski models;
\item c-formulas with $n$ free variables are interpreted by $D^n$-indexed families of 
regular open subsets of $X$, using the usual boolean algebra of regular open sets.
\item $\wn$ is interpreted by taking the interior of the closure, and $\oc$ by doing nothing.
\end{itemize}

\subsection{Subset/sheaf-valued models} \label{sheaf model}

A {\it subset/sheaf-valued structure} interpreting the language of QHC is described as follows.

\begin{itemize}
\item We fix a topological space $X$ and a set $D$;
\item c-formulas with $n$ free variables are interpreted by $D^n$-indexed families of 
arbitrary subsets of $X$, as in the usual Leibniz--Euler--Venn models;
\item i-formulas with $n$ free variables are interpreted by $D^n$-indexed families of 
sheaves of sets on $X$, using the usual operations $\sqcup$, $\x$, $\Hom$, $\bigsqcup_{d\in D}$, 
$\prod_{d\in D}$ on sheaves;
$\Turnstile\Phi$ means that $\Phi$ is interpreted by a family where each sheaf has a global section.
\item $|\oc F|$ is the sheaf given by the inclusion $\Int|F|\to X$;
\item $|\wn\Phi|$ is the set of all $x\in X$ such that the stalk $|\Phi|_x$ of the sheaf $|\Phi|$ is non-empty.
\end{itemize}

In a similar way we describe a subset/{\it pre\/}sheaf-valued structure.
These have somewhat different interpretations of the intuitionistic $\lor$ and $\exists$, 
as $\sigma\F\sqcup\sigma\G\ne\sigma(\F\sqcup\G)$ in general, where $\sigma\F$ denotes 
the presheaf of sections of the sheaf $\F$.

\begin{theorem} \label{global1} Subset/sheaf- and subset/presheaf-valued structures are models of QHC.
\end{theorem}

The results of \S\ref{syntactic interpretations} and \S\ref{principles section} were obtained using 
only the previous three classes of models.
So was the following result: QHC is not complete with respect to either of these 3 classes.

\subsection{Dense image models}

A class of models of QHC with respect to which it is complete was discovered not so long ago by A. Onoprienko
\cite{On1}, \cite{On2}.
Upon seeing her models of the propositional fragment HC \cite{On1}*{\S4}, the author described the following 
new class of {\it topological} models of QHC.

\begin{itemize}
\item Fix a topological space $X$, sets $D$ and $S$ and a map $f\:S\to X$ with dense image;
\item c-formulas with $n$ free variables are interpreted by $D^n$-indexed families of 
arbitrary subsets of $S$, as in the usual Leibniz--Euler--Venn models;
\item i-formulas with $n$ free variables are interpreted by $D^n$-indexed families of
open subsets of $X$, as in the usual Stone--Tang--Tarski models;
\item $\wn$ is interpreted by taking the preimage under $f$;
\item $\oc$ is interpreted by associating to a $T\subset S$ the interior of $X\but f(S\but T)$.
\end{itemize}

The author has asked A. Onoprienko (in 2018) if QHC is complete with respect to this class of topological models.
Recently she answered this question affirmatively in the propositional case, already for those models where $f$ is injective
\cite{On1}*{\S8}, \cite{On3}.

\section{Discussion}

A discussion of related work by other authors can be found in \cite{M1}*{\S\ref{g1:related work}}.

Let us briefly comment on the significance of QHC for elementary mathematics.

1) {\it Algebra.} Subset/sheaf-valued models of QHC (see \S\ref{sheaf model}) arise naturally from consideration of operations on algebraic equations 
\cite{M2}*{\S\ref{g2:equations}}.

2) {\it Geometry.} The logical structure of the first 4 books of Euclid's {\it Elements} is reasonably well described by 
a certain theory over QHC, whose purely classical fragment coincides with Tarski's weak (finitely axiomatized) 
elementary geometry, and purely intuitionistic fragment with Beeson's ``intuitionistic Tarski 
geometry'' \cite{M3}.

3) {\it Arithmetic.} A work in progress by the author introduces a joint logic of propositions, problems and solutions, 
as well as its predicate version QHC$_\lambda$.
It contains QHC as well as S4pr (see \S\ref{proof interpretation}) and the typed $\lambda$-calculus that is related to QH 
by the Curry--Howard correspondence.
If we omit quantification over solutions (and hence the operator $\wn$ of QHC) from QHC$_\lambda$,
the resulting fragment QHC$_\lambda^-$ is expected to admit an interesting arithmetical interpretation:
problems and their solutions are interpreted by arithmetical sets and by their elements as in Kleene's
realizability; and the problem $\oc P$ of proving a proposition $P$ is interpreted by the set of G\"odel numbers
of all formal proofs of $P$.

\bigskip
{\bf Acknowledgments.} 
I would like to thank A.~ Bauer, L.~ Beklemishev, M.~ Bezem, G.~ Dowek, M.~ Jibladze, A.~ Onoprienko, 
A.~ L.~ Semyonov, D.~ Shamkanov, V.~ Shehtman and A.~ Shen' for valuable discussions and useful comments and
the two referees of Doklady for reasonable remarks.
The present note owes its appearance to the requests from A.~ L.~ Semenov and reminders 
from L. Beklemishev and S. Kuznetsov.

%{\bf Disclaimer.}
%I oppose all wars, including those wars that are initiated by governments at the time when 
%they directly or indirectly support my research. The latter type of wars include all wars 
%waged by the Russian state since the Second Chechen war till the present day (October 2, 2022)
%as well as the USA-led invasions of Afghanistan and Iraq.

\end{document}